\documentclass[10pt]{amsart}
\usepackage{color}
\usepackage{soul}
\usepackage{cancel}

\usepackage{amssymb}
\usepackage{ulem}

\definecolor{c20}{rgb}{0.,0.7,0.}
\definecolor{c30}{rgb}{0.,0.,1.}
\definecolor{c40}{rgb}{1,0.1,0.7}
\definecolor{c50}{rgb}{1,0,0}
\definecolor{c60}{rgb}{1,0.9,0.1}

\newcommand{\abs}[1]{\left\lvert #1 \right\rvert}

\newcommand{\E}[1]{\mathbb{E}\left\{#1\right\}}

\newcommand{\pk}[1]{\mathbb{P} \left\{ #1 \right\} }

\newcommand{\R}{\mathbb{R}}

\newcommand{\N}{\mathbb{N}}

\newcommand{\BQN}{\begin{eqnarray}}
\newcommand{\EQN}{\end{eqnarray}}
\newcommand{\BQNY}{\begin{eqnarray*}}
\newcommand{\EQNY}{\end{eqnarray*}}

\newcommand{\BS}{\begin{sat}}
\newcommand{\ES}{\end{sat}}
\newcommand{\BT}{\begin{theo}}
\newcommand{\ET}{\end{theo}}
\newcommand{\BK}{\begin{korr}}
\newcommand{\EK}{\end{korr}}

\newcommand{\BD}{\begin{de}}
\newcommand{\ED}{\end{de}}

\newcommand{\BIT}{\begin{itemize}}
\newcommand{\EIT}{\end{itemize}}
\newcommand{\BDI}{\begin{description}}
\newcommand{\EDI}{\end{description}}

\newcommand{\BRM}{\begin{remarks}}
\newcommand{\ERM}{\end{remarks}}

\newcommand{\BEL}{\begin{lem}}
\newcommand{\EEL}{\end{lem}}

\newtheorem{theo}{Theorem}[section]
\newtheorem{sat}[theo]{Proposition}
\newtheorem{de}[theo]{Definition}
\newtheorem{lem}[theo]{Lemma}

\newtheorem{example}[theo]{Example}
\newtheorem{korr}[theo]{Corollary}
\newtheorem{remark}[theo]{Remark}
\newtheorem{remarks}[theo]{Remarks}
\newtheorem{prop}[theo]{Proposition}

\newcommand{\prooflem}[1]{\textsc{\bf Proof of Lemma} \ref{#1}:}

\newcommand{\COM}[1]{}

\allowdisplaybreaks[1]

\newcommand{\QED}{\hfill $\Box$}

\def\vp{\varepsilon}

\def\IF{\infty}

\def\Cov{\mathrm{Cov}}

\date{}

\def\oo{(1+o(1))}

\def\LT{\left}
\def\RT{\right}
\def\H{\mathcal{H}}

\def\ooo{(1+o(1))}

\def\vn{\varepsilon}

\def\NN{\mathcal{N}}

\def\ovX{\overline X_H}

\def\bnH0{b_n^{H_0/\beta}}

\def\wM{\widetilde{M}}
\def\whM{\widehat{M}}
\def\ws{\widetilde{s}}
\def\wt{\widetilde{t}}
\def\wc{\widetilde{c}}
\def\hc{\hat{c}}
\def\sqln{\sqrt{2\log n}}

\def\Alogn{2A^2\log n}

\usepackage{soul}
\usepackage{color}

\newcommand{\Ji}[1]{{\textcolor{red}{#1}}}

\def\Ji#1{#1}

\begin{document}

\title[
Dependent Gaussian model] 
{On the maxima of suprema of dependent  Gaussian models}

 \author{Lanpeng Ji}
\address{Lanpeng Ji, School of Mathematics, University of Leeds, Woodhouse Lane, Leeds LS2 9JT, United Kingdom
}
\email{l.ji@leeds.ac.uk}

\author{Xiaofan Peng}
\address{Xiaofan Peng, School of Mathematical Sciences, University of Electronic Science and Technology of China, Chengdu 611731, China}
\email{xfpeng@uestc.edu.cn}

\bigskip

\date{\today}
 \maketitle

{\bf Abstract:} In this paper, we study the asymptotic distribution of the maxima of suprema of dependent Gaussian processes with trend. For different scales of the time horizon we obtain different normalizing functions for the convergence of the maxima. The obtained results not only 
have  potential applications in estimating the delay of certain Gaussian fork-join queueing systems but also provide interesting insights to the extreme value theory for triangular arrays of random variables with row-wise dependence.
\medskip

\medskip


\section{Introduction}

Let  $\{X_i(t), t\ge0\}, i=1,2,\ldots,$ be independent copies of a centered self-similar Gaussian process with almost surely (a.s.) continuous sample paths,   self-similarity index $H\in(0,1)$ and variance function $t^{2H}$, and let  $\{X(t), t\ge0\}$ be another independent  centered self-similar Gaussian processes with a.s. continuous sample paths,  self-similarity index $H_0\in(0,1)$ and variance function $t^{2H_0}$. 
We define, for positive constants $\sigma,\sigma_0, c_i, i\ge1$, $\beta>\max(H,H_0)$, and a deterministic function $T_n>0$,
\BQN\label{eq:Mn}
M_n:=  \underset{i\le n}{\max} \ \sup_{t\in[0,T_n]}(\sigma X_i (t)  +\sigma_0 X(t) -c_i t^\beta),\ \ \ \ n\ge 1.
\EQN
This paper is concerned with asymptotic distributional properties of $M_n$ as $n\to\IF.$ More precisely, we aim to establish limit theorems for $\nu_n^{-1}(M_n-\mu_n)$, as $n\to\IF$, for some suitably chosen normalizing functions $\nu_n, \mu_n, n\ge1.$ This work is a continuation of the recent work done in \cite{JP22}, where the case $T_n=\IF$ was discussed. Note that the general $k$th order statistics were discussed in \cite{JP22}, but to ease the complication 
we shall only focus on the maxima in this paper; one can check that some results for the $k$th order statistics can also be established. As in \cite{JP22}, without loss of generality, we assume $\sigma=1$. Throughout the rest of the paper, we also assume that $\lim_{n\to\IF} T_n \in[0,\IF]$ exists and investigate how different scales of $T_n$ influence the normalizing functions $\nu_n, \mu_n, n\ge1$ in the limit theorems for  the maxima defined in \eqref{eq:Mn}. 

The motivation for the study of distributional properties of $M_n$ stems from a recent contribution \cite{MSJVZ21} on a Brownian fork-join queueing system, which  is a special model of \eqref{eq:Mn} with all the Gaussian processes involved being Brownian motions, $\beta=1$, $c=c_i,i\ge1$ and $T_n=\IF$ (hereafter called {\it Browian model with linear drift}).  In their context, $M_n$ (with $T_n=\IF$)
models the maximum of stady-state queue lengths (or delay)  in
a fork-join network of $n$ statistically identical queues  which are driven by a common Brownian motion perturbed arrival process and  independent  Brownian motion perturbed service processes, respectively. The  theoretical limit result obtained therein is the key to developing structural insights into the dimensioning of assembly systems; interested readers are referred to \cite{MSJVZ21} for more details on this application. More recently, the tail asymptotics for the delay in such a Brownian fork-join queueing system were studied in \cite{SVZ22}.
As discussed in \cite{DM03, DR02,  Man07} and references therein, for a fluid queueing model it is of great interest to consider general Gaussian processes with a non-linear trend and study the distributional properties of the transient queue length (i.e., the supremum is taken over a finite time interval instead of $\R_+$). 
Analogously,
the maximum $M_n$ in \eqref{eq:Mn} can be seen as the maximum of $n$ transient queue lengths  (or delay over a finite time horizon) in a general Gaussian fork-join queueing system, with the time horizon $T_n$  possibly dependent on $n$. In this paper, we 
consider different scales of $T_n$. This discussion may be intersting from an application point of view, for instance, the system users may be interested in estimating the  delay over any short-time  or long-time horizon. 
Note that as in \cite{JP22} the study in this paper also provides complementary results to the extreme value theory for multivairate Gaussian models in random environment. 



Clearly, the study of the maxima $M_n, n\ge 1$ is relevent to the extreme value theory for triangular arrays of random variables with row-wise dependence.  Define $\{Y_{kn}, k\le n, n\ge1\}$ to be a triangular array of random variables and $N_n=\max_{k\le n} Y_{kn}$ to be the row-wise maximum. The extreme value theory for the triangular array  $\{Y_{kn}, k\le n, n\ge1\}$ is concerned with the convergence of the row-wise maxima $N_n, n\ge1$ under a linear normalization. If $Y_{kn}, k\le n$ is stationary for any fixed  $n$, then we call $N_n, n\ge1$ {\it homogeneous} maxima, otherwise, we call it {\it inhomogeneous} maxima.
Current literature on extreme value theory for triangular arrays has been focused on homogeneous maxima, and particularly, the maxima for row-wise independent and indentically distributed triangular arrays (i.e., $Y_{kn}, k\le n$ being independent and identically distributed (IID)); 
see, e.g.,  \cite{ACH97, FH03, PM21}. Some conditions guaranteeing the convergence of normalized maxima to some limit are given by \cite{FH03} under some differentiation conditions. 
Particularly, the maxima of  row-wise independent  Poisson-distributed and related triangular arrays are discussed in \cite{ACH97}, and the maxima for some  row-wise independent  Weibull-(truncated) regular variation mixture distributed triangular arrays are discussed in \cite{PM21}. Somehow surprisingly, except for the triangular arrays of normal random variables (e.g., \cite{HHR96ab}), there are very few papers dealing with the homogeneous maxima with row-wise stationary triangular arrays. 
The only result on this topic that we could find is the recent one obtained in \cite{DEN16} where a Gumbel limit theorem is obtained for normalized maxima under some general conditions (see Theorem 2.1 therein).  
It turns out that  there exists no   theory for general (in)homogeneous maxima which covers the convergence of $M_n, n\ge1$ under some normalization that is interested in this paper.
 In what follows, if $c=c_i, i\ge1$, the maxima $M_n, n\ge 1$ is called {\it homogeneous}, and otherwise, called {\it inhomogeneous}.
 We obtain convergence results for suitably normalized $M_n,n\ge 1$  for both homogeneous case and some inhomogeneous case. As we will see, the only possible non-degenerate limit distributions are from the family of Gumbel, Gaussian or a mixture of them. This study  provides some interesting examples, which enriches the extreme value theory for triangular arrays of random variables with row-wise dependence.

\medskip

The rest of the paper is organized as follows: In Section 2 we  present some preliminary results concerning the tail asymptotics of the supremum of  a class of self-similar Gaussian processes with trend over a threshold-dependent time horizon. The main results on the homogeneous maxima $M_n,n\ge 1$ are given in Section 3.  Section 4 discusses some   inhomogeneous maxima. All the proofs are presented in Section 5.

\section{Prelimanaries} 
In this section, we mainly discuss the tail asymptotics of the supremum of a self-similar Gaussian process with trend over a threshold-dependent time interval. This study is useful for the construction of normalizing functions for the maxima, and is also of independent interest.  The results presented in Proposition \ref{propTn} below generalize some of the existing results obtained in \cite{BPP09}, see also \cite{LHJ15}.

Let $\{X_H(t), t\ge0\}$ be a centered self-similar Gaussian process with a.s. continuous sample paths, self-similarity index $H\in(0,1),$ variance function $t^{2H}$. We assume  a {\it local stationarity}
of the standardized  Gaussian process $\ovX(t):=X_H(t)/t^H, t> 0$ in a
neighbourhood of the point $t=1$, i.e.,
\BQN\label{eq:locstaXH}
\lim_{s,t\to 1}\frac{ \E{(\ovX(s)-\ovX(t))^2}}{K^2 (\abs{s-t})}=1
\EQN
holds for some positive function $K(\cdot)$ which is regularly varying at 0 with index $\alpha/2\in(0,1)$. Condition \eqref{eq:locstaXH} is a common assumption in the literature; see, e.g., \cite{dieker2005extremes} and \cite{HP99}.
It is worth noting that the assumption \eqref{eq:locstaXH} is slightly general than the S2 
 in \cite{DT20} where a decent discussion on properties and examples of self-similar Gaussian processes is given.  
Note that the local stationarity at $t=1$ and the self-similarity of the random process imply the local stationarity at any point $t=r>0$, i.e.,
\BQN\label{eq:locstaXHr}
\lim_{s,t\to r}\frac{ \E{(\ovX(s)-\ovX(t))^2}}{K^2 (\abs{s-t})}=r^{-\alpha}.
\EQN

For a threshold-dependent time horizon $T_u$ (to be specified below) and constants $c>0, \beta>H$, we shall derive asymptotics for
\BQNY
\psi_{T_u}(u):=\pk{\sup_{t\in[0,T_u]} X_H(t) -ct^\beta >u   },\  \ \ \ u\to\IF.
\EQNY
Throughout this paper, for two positive functions $f, h$ and some $u_0\in[-\IF, \IF]$, 
write $h(u)\sim  f(u)$ or $h(u)=f(u)\oo$ if $ \lim_{u \to u_0} f(u) /h(u)  = 1 $,  write $ f(u) = o(h(u)) $ if $ \lim_{u \to u_0} {f(u)}/{h(u)} = 0$, and  write $ f(u) =O (h(u)) $ if $ \lim_{u \to u_0} {f(u)}/{h(u)} \in(0,\IF)$.
Further, we denote by $ \overset{\leftarrow}{K}(\cdot)$  the asymptotic inverse of  $K(\cdot)$, and thus 
 $$
 \overset{\leftarrow}{K}(K(t))=K( \overset{\leftarrow}{K}(t))\oo=t\oo,\ \ \ t\downarrow0.
  $$
It follows that  $ \overset{\leftarrow}{K}(\cdot)$ is regularly varying at 0 with index ${2}/{\alpha}$; see, e.g., \cite{EKM97}.

We shall consider the following scenarios for the threshold-dependent time horizon $T_u$:
\begin{itemize}
\item[\bf{D1:}] 
 $\lim_{u\to\IF} T_u/u^{1/\beta}=0$;
\item[\bf{D2:}] 
$\lim_{u\to\IF} T_u/u^{1/\beta}=s_0\in(0,t_0)$;
\item[\bf{D3:}] 
$\lim_{u\to\IF}\frac{T_u-t_0u^{1/\beta}}{A^{1/2} B^{-1/2}u^{H/\beta+1/\beta-1}}=x\in (-\IF, \IF]$.
\end{itemize}
Here in the above 
$$
t_0=\LT(\frac{H}{c(\beta-H)}\RT)^{1/\beta},\ \ \ \ 
$$
and
\BQN\label{eq:AB}
 A=\frac{t_0^H}{1+ct_0^\beta}=\frac{\beta-H}{\beta}\LT(\frac{H}{c(\beta-H)}\RT)^{H/\beta},\ \ \ \ B=\LT(\frac{H}{c(\beta-H)}\RT)^{-\frac{H+2}{\beta}} H \beta.
\EQN

Below, by $\{B_{\alpha/2}(t), t\ge 0\}$ we denote a standard fractional Brownian motion
 (sfBm)
with Hurst index $\alpha/2\in(0,1)$, and 
\BQNY
\Cov(B_{\alpha/2}(t),B_{\alpha/2}(s))=\frac{1}{2}(t^{\alpha}+s^{\alpha}-\mid t-s\mid^{\alpha}),\ \ \ t,s\ge0.
\EQNY
The well known Pickands constant $\H_{\alpha}$ and Piterbarg constant $\mathcal{P}_\alpha^{d}$ in the Gaussian theory is defined, respectively, by
\BQNY
&&\H_{\alpha}=\lim_{T\to\infty}\frac{1}{T}
\E{\exp \left(\sup_{t\in[0,T]}(\sqrt{2}B_{\alpha/2}(t)-t^\alpha)\right)} \in (0,\IF).\\
&&\mathcal{P}_\alpha^{d}=\lim_{T\to\infty}
\E{\exp \left(\sup_{t\in[0,T]}(\sqrt{2}B_{\alpha/2}(t)-(1+d)t^\alpha)\right)} \in (0,\IF),\ \ \ \ d>0.
\EQNY
We refer to  \cite{BJ22,DH20, demiro2003simulation,  DikerY, Pit96}
for basic properties of {the} Pickands, Piterbarg and related constants. In particular, it has been shown that $\mathcal{H}_{1}=1$ and $\mathcal{P}_1^{d}=1+1/d,$ $d>0.$


\begin{prop} \label{propTn}
Let $\{X_H(t),t\ge 0\}$ be the self-similar Gaussian process defined as above with \eqref{eq:locstaXH} satisfied, and assume $c>0,$ $\beta>H$.

\begin{itemize}
\item[(i).]
Further assume that the following limit exists and
\BQNY
\lim_{t\downarrow 0}\frac{t}{K^2(t)}=:Q\in[0,\IF].
\EQNY
Then,  under scenarios D1 and D2 (i.e., $\lim_{u \rightarrow \infty} {T_u}/{u^{1/\beta}} = s_0 \in[0, t_0)$), we have
\BQNY 
\psi_{T_u}(u)=\mathcal{D}_{c_0}\left(\frac{u+cT_u^\beta}{T^H_u}\right)\left(\frac{u+cT_u^\beta}{T^H_u}\right)^{-1}\exp \left(-\frac{(u+cT_u^\beta)^2}{2T^{2H}_u}\right)\ooo, \quad u\to \IF,
\EQNY
where, for $y>0,$
\BQN\label{eq:DH}
\mathcal{D}_{c_0}(y)=\left\{
            \begin{array}{ll}
  \frac{\mathcal{H}_{\alpha}}{2^{1/\alpha}\sqrt{2\pi}(H-c_0\beta)} y^{-2}\left( \overset{\leftarrow}K\LT(1/y\RT)\right)^{-1} , & \hbox{if } Q=0 ,\\
\frac{1}{\sqrt{2\pi}}\  \mathcal{P}^{2(H-c_0\beta)Q}_\alpha, & \hbox{if } Q\in(0,\IF),\\
\frac{1}{\sqrt{2\pi}}&  \hbox{if } Q=\IF,
              \end{array}
            \right.\ \ \ \text{with}\ c_0=\frac{cs_0^\beta}{1+cs_0^\beta}.
\EQN
\item[(ii).] Under scenario D3  (i.e., $\lim_{u\to\IF}\frac{T_u-t_0u^{1/\beta}}{A^{1/2} B^{-1/2}u^{H/\beta+1/\beta-1}}=x \in (-\IF, \IF]$), we have
\BQNY 
\psi_{T_u}(u)=\psi_{\IF}(u)\Phi(x)\ooo, \quad u\to \IF,
\EQNY
where $\Phi(x)$ is the standard normal distribution function and
\BQNY
\psi_{\IF}(u):=\pk{\sup_{t\ge0} X_H(t) -ct^\beta >u   }=R(u)\exp\LT(-\frac{u^{2\LT(1-\frac{H}{\beta}\RT)}}{2A^2}\RT)\oo,\ \ u\to\IF,\label{eq:HP}
\EQNY
where (with $A,B$ given  in \eqref{eq:AB}) 
\BQN\label{Ru}
R(u)=\frac{A^{\frac{3}{2}-\frac{2}{\alpha}}  \H_\alpha}{2^{ \frac{1}{\alpha}} B^{\frac{1}{2}} t_0}\frac{ u^{  \frac{2H}{\beta}-2}}{ \overset{\leftarrow}K (u^{\frac{H}{\beta}-1}) }, \ \  u>0.
\EQN

\end{itemize}

\end{prop}

\begin{remark}
 If $\{X_H(t), t\ge0\}$ is the sfBm with Hurst index $H\in(0,1)$ and $\beta=1$, then we can check that $K(t)=t^H$, and the results in Proposition \ref{propTn} reduce to that given in (2.2)-(2.5) in \cite{LHJ15} (here the explicit formula for $\mathcal{P}_1^{d}$ given above should be used). 


\end{remark}

\section{ Homogeneous maxima}

We shall consider asymptotic distributional properties of the homogeneous  maxima 
$$
M_n=  \underset{i\le n}{\max} \ \sup_{t\in[0,T_n]}(X_i (t)  +\sigma_0 X(t) -c t^\beta), \ \ \ \ n\ge1,
$$
where we assume $\sigma=1$ and $c=c_i, i\ge1.$
Precisely, we aim to establish limit theorems for $\nu_n^{-1}(M_n-\mu_n)$, as $n\to\IF$, for some suitably chosen normalizing functions $\nu_n, \mu_n, n\ge1.$
Motivited by the scenarios D1-D3 treated in Section 2, we consider the following scenarios for $T_n$:

\begin{itemize}
\item[\bf{S1.}] ({\it Super-short time horizon}) $\lim_{n\to\IF}T_n^{H} \sqrt{2\log n}= \kappa_0 \in [0,\IF);$
\item[\bf{S2.}]  ({\it Short time horizon})  $\lim_{n\to\IF}T_n^{H} \sqrt{2\log n}=\IF$ and $\lim_{n\to\IF} \frac{T_n^{1-H/\beta}}{ (\sqrt{2\log n})^{1/\beta}}=0$;
\item[\bf{S3.}]  ({\it Intermediate time horizon})  $\lim_{n\to\IF} \frac{T_n^{1-H/\beta}}{ (\sqrt{2\log n})^{1/\beta}}=\ws_0\in\LT(0,\wt_0\RT)$, with $\wt_0=\LT(\frac{t_0^\beta}{1+ct_0^\beta}\RT)^{1/\beta}=\LT(\frac{H}{c\beta}\RT)^{1/\beta}$;
\item[\bf{S4.}]  ({\it Long time horizon}) $\lim_{n\to\IF} \frac{T_n^{1-H/\beta}}{ (\sqrt{2\log n})^{1/\beta}}=\wt_0$ and $\lim_{n\to\IF}\frac{T_n- (\wt_0^\beta \sqrt{2\log n})^{1/(\beta-H)}}{A^{1/2} B^{-1/2} (2A^2\log n)^{(H+1-\beta)/(2(\beta-H))} }=x_0\in \R$;
\item[\bf{S5.}] ({\it Super-long time horizon}) $\lim_{n\to\IF} \frac{T_n^{1-H/\beta}}{ (\sqrt{2\log n})^{1/\beta}} =\ws_0>\wt_0$.
\end{itemize}

To determine the correct normalizing functions, we will discuss the following maxima of the corresponding IID sequence
$$
 \wM_n :=\underset{i\le n}{\max} \ \sup_{t\in[0,T_n]}(X_i (t)   -c t^\beta),\ \ \ \ n\ge 1.
$$
It is known that normalizing functions for the maxima of a sequence of IID random variables can be retrieved from the tail asymptptics of the random variable; see, e.g., \cite{JP22}. Motivited by the tail asymptotics discussed in Proposition \ref{propTn}, we introduce the following normailizing functions.

Under scenarios S2 and S3:
\BQN\label{eq:cndn}
b_n
&:=&T_n^H\LT(\sqln + \frac{1}{\sqln} \log \LT(f_n\LT(T_n^H\sqln -c T_n^\beta\RT)\RT)\RT) -c T_n^\beta,\  n\in \N,\\
a_n& :=& 
\frac{ T_n ^{H} }{\sqrt{2  \log n}} , \ n\in \N.   \nonumber
\EQN
where, recalling the function $\mathcal D_{c_0}(\cdot)$ defined in \eqref{eq:DH},
\BQN \label{eq:tauC}
 f_n(w):=\mathcal{D}_{c \ws_0^\beta}\left(\frac{w+cT_n^\beta}{T^H_n}\right)\left(\frac{w+cT_n^\beta}{T^H_n}\right)^{-1},  \ \ w>0,  
\EQN
Under scenarios S4 and S5 (for S5, we set $x_0=\IF$, i.e., use $d_n(\IF)$):
\BQN\label{eq:bnanS3}
d_n(x_0)& :=& 
 (2A^2 \log n)^{1/\tau} \LT(1 +\frac{1}{\tau  \log n}  \log(R((2A^2\log n)^{1/\tau})\Phi(x_0))  \RT)
,\  n\in \N,\\
e_n& :=&
\frac{ (2A^2 \log n)^{1/\tau}}{\tau \log n}, \ n\in \N,  \nonumber
\EQN
where $ \tau = 2(1-H/\beta),$ and  $R(u)$ is defined in \eqref{Ru}.

The following proposition shows that these normalizing functions are the required ones for $ \wM_n, n\ge 1$. Hereafter, $\overset{d}\to$ denotes the convergence in distribution, and $\Lambda$ denotes a standard Gumbel random varaiable, i.e., $\pk{\Lambda\le x} = \exp(-e^{-x}), x\in\R.$
\begin{prop} \label{prop:wMbS}
We have
\begin{itemize}
\item[(i).] Under scenarios S2 and S3,
\BQNY 
a_n^{-1}(\wM_n  -b_n ) \overset{d}\to \Lambda, \ \ \ \  n\to\IF.
\EQNY
\item[(ii).] Under scenarios S4 and S5,
\BQNY 
e_n^{-1}(\wM_n  -d_n(x_0) ) \overset{d}\to \Lambda,\ \ \ \  n\to\IF.
\EQNY

\end{itemize}

\end{prop}

Below is our pricinpal result on the aymptotical distribution for the homogenous maxima $M_n, n\ge1$ under suitable normalization. In what follows, we denote  $\NN$ to be a standard normal random variable which is independent of $\Lambda$.

\BT \label{Thm:main1}

For the homogeneous maxima $M_n, n\ge 1$ defined in \eqref{eq:Mn} with $\sigma=1$ and $c=c_i, i\ge 1$, we have

\begin{itemize}
\item[(a).] Under scenario S1,
\BQNY
M_n \overset{d}\to \kappa_0, \ \ \ \  n\to\IF.
\EQNY

\item[(b).] Under scenarios S2 and S3,

\begin{itemize}
\item[(b.i)] If further $\lim_{n\to\IF} \frac{T_n^{H-H_0}}{ \sqrt{2\log n}}=0$, then
\BQNY
\sigma_0^{-1}T_n^{-H_0} (M_n - b_n) \overset{d}\to \NN, \ \ \ \  n\to\IF.
\EQNY

\item[(b.ii)] If further $\lim_{n\to\IF} \frac{T_n^{H-H_0}}{ \sqrt{2\log n}}=\IF$, then
\BQNY
a_n^{-1} (M_n - b_n) \overset{d}\to \Lambda, \ \ \ \  n\to\IF.
\EQNY

\item[(b.iii)] If further $\lim_{n\to\IF} \frac{T_n^{H-H_0}}{ \sqrt{2\log n}}=q_0\in(0,\IF)$, then
\BQNY
a_n^{-1} (M_n - b_n) \overset{d}\to \Lambda +  \frac{\sigma_0}{q_0} \NN, \ \ \ \  n\to\IF.
\EQNY
\end{itemize}

\item[(c).] Under scenarios S4 and S5,

\begin{itemize}
\item[(c.i)] If further $2H-H_0<\beta$,
then
\BQNY
\sigma_0^{-1}T_n^{-H_0} (M_n - d_n(x_0)) \overset{d}\to \NN, \ \ n\to\IF.
\EQNY

\item[(c.ii)] If further $2H-H_0>\beta$, then
\BQNY
e_n^{-1} (M_n - d_n(x_0)) \overset{d}\to \Lambda, \ \ n\to\IF.
\EQNY

\item[(c.iii)] If further $2H-H_0=\beta$, then
\BQNY
e_n^{-1} (M_n - d_n(x_0)) \overset{d}\to \Lambda + \frac{\sigma_0 c\beta  }{H} \NN, \ \ n\to\IF.
\EQNY
\end{itemize}

\end{itemize}

\ET

\begin{remarks}
(a). We remark that the conditions of (b.i)-(b.iii) can be more specific under scenario S3. In fact, under S3 we know that
\BQNY
T_n \sim  \LT(\ws_0^\beta \sqrt{2\log n}\RT )^{ \frac1{\beta-H}}, \ \ \ \ n\to\IF,
\EQNY
and thus the case-specific conditions of (b.i)-(b.iii) can be simplified to $2H-H_0<\beta$,  $2H-H_0>\beta$ and $2H-H_0=\beta$, respectively, and further $q_0 = \ws_0^\beta. $

(b). 
It turns out that there is an interesting smooth transition from $a_n, b_n$, respectively, to $e_n, d_n(x_0)$ in the sense that, if $T_n$ satisfies S4 then
\BQNY
b_n\sim T_n^H \sqln   -cT_n^\beta \sim (2A^2 \log n)^{1/\tau}\sim d_n(x_0), \ \ \text{and} \ \ a_n \sim e_n, \ \ \ \ \ n \to\IF.
\EQNY
In fact, 
we can rewrite $T_n$ as
\BQN\nonumber
T_n&=&\LT(\wt_0^\beta\sqrt{2\log n}\RT)^{\frac{1}{\beta-H}}+A^{1/2} B^{-1/2}x_0  (2A^2\log n)^{(H+1-\beta)/(2(\beta-H))}+\vp(n)\\\label{def:Tn-s4}
&=& t_0(\Alogn)^{\frac{1}{2(\beta-H)}} \LT(1+ \frac{A^{1/2} B^{-1/2}x_0}{t_0}(\Alogn)^{-1/2}  + t_0^{-1}\vp(n) (\Alogn)^{\frac{-1}{2(\beta-H)}}\RT),
\EQN
where $\vp(n)=o\LT( (\log n)^{(H+1-\beta)/(2(\beta-H))}\RT)$. Thus,
\BQNY
&&\quad T_n^{H}\sqln-cT_n^\beta - (\Alogn)^{1/\tau}\\
&&= A^{-1} t_0^H(\Alogn)^{\frac{\beta}{2(\beta-H)}} \LT( 1+ H t_0^{-1}A^{1/2} B^{-1/2}x_0 (\Alogn)^{-1/2}  + H t_0^{-1}\vp(n) (\Alogn)^{\frac{-1}{2(1-H)}} +O((\Alogn)^{-1}) \RT)\\
&&\quad - c t_0 (\Alogn)^{\frac{\beta}{2(\beta-H)}} \LT( 1+ \beta t_0^{-1}A^{1/2} B^{-1/2}x_0 (\Alogn)^{-1/2}  + \beta t_0^{-1}\vp(n) (\Alogn)^{\frac{-1}{2(1-H)}} +O((\Alogn)^{-1}) \RT)\\
&&\ \ \ -(\Alogn)^{1/\tau} \\
&&
=o\LT((\Alogn)^{1/\tau}\RT), \ \ \ \  n\to\IF,
\EQNY
where in the second equality we have used $A^{-1}t_0^{H}=1+ct_0$. 
Similarly, some elementary calculations show that 
$$
\lim_{n\to\IF}{a_n}/{e_n}=1.
$$  This is an  interesting observation which reveals a smooth change of the normalising functions.
However, it looks that under scenario S4 
 the $b_n$ is not the correct normalising function but the $d_n(x_0)$ is. 

\end{remarks}

\COM{

\begin{remarks}
(a). An example of $T_n$ that satisfies the assumptions of  scenario 1) which allows for three different further scenarios (i)-(iii) is given by
$$
T_n = (\sqrt{2 \log n})^{\frac{1}{1-H} -\delta},\ \ \text{with} \ \delta\in\LT(0, \frac{1}{1-H}\RT).
$$

(c...)

\end{remarks}
} 

\begin{example} \label{example:1}
This example illustrates  Theorem \ref{Thm:main1} for $T_n=(\lambda \sqln)^\gamma$, with $\gamma\in\R$ and $\lambda>0$. We check how the asymptotic distribution of a normalized $M_n$ will change according to different values of $\gamma$. There can be a lot of cases to be considered, but we choose to work with one representative case where $\beta>H_0>2H$. In this case, we obtain different results according to where the value of $\gamma$ falls in the following intervals:
$$
-\IF<-\frac1{H}<-\frac{1}{H_0-H}<\frac{1}{\beta-H}<\IF.
$$ 
Precisely, a direct appliation of Theorem \ref{Thm:main1} yeilds the following convergence results, as $n\to\IF$:
\begin{itemize}
\item[(1).] If $\gamma\in\LT(-\IF, -\frac1{H}\RT)$, then $M_n \overset{d}\to 0$.
\item[(2).] If $\gamma=-\frac1{H}$, then $M_n \overset{d}\to \lambda^{-1}$.
\item[(3).] If $\gamma\in\LT( -\frac1{H}, -\frac{1}{H_0-H}\RT)$, then $a_n^{-1} (M_n - b_n) \overset{d}\to \Lambda$.
\item[(4).] If $\gamma=-\frac{1}{H_0-H}$, then $a_n^{-1} (M_n - b_n) \overset{d}\to \Lambda+  \sigma_0 \lambda^{-1} \NN.$
\item[(5).] If $\gamma\in\LT( -\frac{1}{H_0-H}, \frac{1}{\beta-H}\RT)$, then $\sigma_0^{-1}T_n^{-H_0}(M_n - b_n) \overset{d}\to \NN$.
\item[(6).] If $\gamma= \frac{1}{\beta-H}$ and $\lambda\in(0,\wt_0^\beta)$, then $\sigma_0^{-1}T_n^{-H_0} (M_n - b_n) \overset{d}\to \NN$.
\item[(7).] If $\gamma= \frac{1}{\beta-H}$ and $\lambda=\wt_0^\beta$, then $\sigma_0^{-1}T_n^{-H_0} (M_n - d_n(0)) \overset{d}\to \NN$.
\item[(8).] If $\gamma= \frac{1}{\beta-H}$ and $\lambda\in(\wt_0^\beta,\IF)$, or $\gamma> \frac{1}{\beta-H}$,  then $\sigma_0^{-1}T_n^{-H_0} (M_n - d_n(\IF)) \overset{d}\to \NN$.
\end{itemize}

\end{example}

\section{  Inhomogenous maxima}

We shall consider convergence results for the inhomogenous maxima 
$$
M_n=  \underset{i\le n}{\max} \ \sup_{t\in[0,T_n]}(X_i (t)  +\sigma_0 X(t) -c_i t^\beta),\ \ \ \ n\ge1,
$$
where we assume $\sigma=1.$
For simplicity,  we shall assume that all the $c_i, i\ge1$ take value from a finite set of distinct values, denoted as 
\BQN\label{eq:Sc}
\mathcal S=\{\hc_1, \hc_2,\cdots,\hc_k\},\ \ \text{ where }\ 0<\hc_1(=c)<\hc_2<\cdots<\hc_k<\IF,
\EQN with some $k\ge2.$
Here, for notational convenience, we use $\hc_1=c$ which is helpful when we adopt the normalizing functions defined in the previous sections.
Further, let $m_j=\#\{i\le n: c_i = \hc_j\}$ and  assume that
\BQN\label{eq:pp}
\lim_{n\to\IF} \frac{m_1}{n}=:p_1\in(0,1], \  \ \lim_{n\to\IF} \frac{m_j}{n}=:p_j\in[0,1), \ \ 2\le  j\le k.
\EQN
Obviously, $\sum_{j=1}^k p_j=1.$

Similarly to the homogenous case we shall first discuss the maxima of the corresponding IID sequence
\BQN\label{eq:whM}
 \whM_n :=\underset{i\le n}{\max} \ \sup_{t\in[0,T_n]}(X_i (t)   -c_i t^\beta),\ \ \ \ n\ge 1.
\EQN
We have the following result for the asympototical distribution of suitably normalized maxima $\whM_n, n\ge 1$.

\begin{prop} \label{prop:wMbSc}
Let $\whM_n, n\ge 1$ be defined as in \eqref{eq:whM}, with $c_i, i\ge 1$ satisfying \eqref{eq:Sc} and \eqref{eq:pp}. We have

\begin{itemize}
\item[(a).] Under S1, 
\BQNY
\whM_n  \overset{d}\to \kappa_0,\ \ \ \ n\to\IF.
\EQNY

\item[(b).]  Under S2, 
\BQNY
a_n^{-1}(\whM_n  -b_n ) \overset{d}\to \widehat \Lambda, \ \ \ \  n\to\IF,
\EQNY
where
\BQNY
 \widehat \Lambda = \left\{
\begin{array}{lll}
\Lambda, &  \mbox{if\ } \lim_{n\to\IF} T_n^{\beta-H} \sqrt{2\log n}=0,\\[0.1cm]
\Lambda+ \log\LT(p_1 +\sum_{j=2}^k p_j e^{-(\hc_j-c)q_1}\RT), &  \mbox{if\ } \lim_{n\to\IF} T_n^{\beta-H} \sqrt{2\log n}=q_1\in(0,\IF),\\[0.1cm]
\Lambda+\log p_1, & \mbox{if\ } \lim_{n\to\IF} T_n^{\beta-H} \sqrt{2\log n}=\IF.
\end{array}
\right.
\EQNY

\item[(c).]  Under S3, 
\BQNY
a_n^{-1}(\whM_n  -b_n ) \overset{d}\to \Lambda+\log p_1,\ \ \ \ n\to\IF.
\EQNY
\item[(d).] Under S4 and S5, 
\BQNY
e_n^{-1}(\whM_n  -d_n(x_0) ) \overset{d}\to \Lambda+\log p_1,\ \ \ \ n\to\IF.
\EQNY


\end{itemize}

\end{prop}

\begin{remarks}
(a). It is interesting to observe that under S2, the three possible limits (i.e., three values of $\widehat \Lambda$) are all from the Gumbel family, where the second one depends on all the constants $\hc_j, 1\le j\le k$,  the third one depends only on the proportion of $\hc_1=c$, and the first one is not really affected by the more specific information of the trend functions. 

(b). After some algebric calculations for the normalizing functions, one can check that the result for S5 is consistent with the result of Theorem 3.5 in \cite{JP22}.

\end{remarks}

Below is the main result of this section.

\BT
\label{Thm:main2}
Let $M_n, n\ge 1$ be the inhomogenous maxima defined as in \eqref{eq:Mn}, with $\sigma=1$, $c_i, i\ge 1$ satisfying \eqref{eq:Sc} and \eqref{eq:pp}. We have

\begin{itemize}
\item[(a).] Under S1, 
\BQNY
M_n  \overset{d}\to \kappa_0,\ \ \ \ n\to\IF.
\EQNY

\item[(b).]  Under S2,  
\begin{itemize}
\item[(b.i)] If further $\lim_{n\to\IF} \frac{T_n^{H-H_0}}{ \sqrt{2\log n}}=0$, then
\BQNY
\sigma_0^{-1}T_n^{-H_0} (M_n - b_n) \overset{d}\to \NN, \ \ \ \  n\to\IF.
\EQNY

\item[(b.ii)] If further $\lim_{n\to\IF} \frac{T_n^{H-H_0}}{ \sqrt{2\log n}}=\IF$, then
\BQNY
a_n^{-1} (M_n - b_n) \overset{d}\to  \widehat \Lambda, \ \ \ \  n\to\IF.
\EQNY

\item[(b.iii)] If further $\lim_{n\to\IF} \frac{T_n^{H-H_0}}{ \sqrt{2\log n}}=q_0\in(0,\IF)$, then
\BQNY
a_n^{-1} (M_n - b_n) \overset{d}\to  \widehat \Lambda+  \frac{\sigma_0}{q_0} \NN, \ \ \ \  n\to\IF,
\EQNY
with $\widehat\Lambda$ defined as in Proposition \ref{prop:wMbSc}(b).
\end{itemize}
\COM{\BQNY
 \widetilde \Lambda = \left\{
\begin{array}{lll}
\Lambda, &  \mbox{if\ } \lim_{n\to\IF} T_n^{\beta-H} \sqrt{2\log n}=0,\\[0.1cm]
\Lambda+ \log\LT(p_1 +\sum_{j=2}^k p_j e^{-(\hc_j-c)q_1}\RT), &  \mbox{if\ } \lim_{n\to\IF} T_n^{\beta-H} \sqrt{2\log n}=q_1\in(0,\IF),\\[0.1cm]
\Lambda+\log p_1, & \mbox{if\ } \lim_{n\to\IF} T_n^{\beta-H} \sqrt{2\log n}=\IF.
\end{array}
\right.
\EQNY}

\item[(c).]  Under S3, 

\begin{itemize}
\item[(c.i)] If further $\beta>2H-H_0$, then
\BQNY
\sigma_0^{-1}T_n^{-H_0} (M_n - b_n) \overset{d}\to \NN, \ \ \ \  n\to\IF.
\EQNY

\item[(c.ii)] If further $\beta<2H-H_0$, then
\BQNY
a_n^{-1} (M_n - b_n) \overset{d}\to  \Lambda+\log p_1, \ \ \ \  n\to\IF.
\EQNY

\item[(c.iii)] If further $\beta=2H-H_0$, then
\BQNY
a_n^{-1} (M_n - b_n) \overset{d}\to \Lambda+\log p_1+  \frac{\sigma_0}{\ws_0^\beta} \NN, \ \ \ \  n\to\IF.
\EQNY
\end{itemize}

\item[(d).] Under S4 and S5, 
\begin{itemize}
\item[(d.i)] If further $\beta>2H-H_0$, then
\BQNY
\sigma_0^{-1}T_n^{-H_0} (M_n - d_n(x_0)) \overset{d}\to \NN, \ \ \ \  n\to\IF.
\EQNY

\item[(d.ii)] If further $\beta<2H-H_0$, then
\BQNY
e_n^{-1} (M_n - d_n(x_0)) \overset{d}\to  \Lambda+\log p_1, \ \ \ \  n\to\IF.
\EQNY

\item[(d.iii)] If further $\beta=2H-H_0$, then
\BQNY
e_n^{-1} (M_n - d_n(x_0)) \overset{d}\to \Lambda+\log p_1+  \frac{\sigma_0c\beta}{H} \NN, \ \ \ \  n\to\IF.
\EQNY
\end{itemize}


\end{itemize}

\ET

\begin{remark}
In Theorem \ref{Thm:main2}(b), we introduce mixture conditions according to the possible limit values of 
$\lim_{n\to\IF} \frac{T_n^{H-H_0}}{ \sqrt{2\log n}}$ and $\lim_{n\to\IF} T_n^{\beta-H} \sqrt{2\log n}$ (see also Proposition \ref{prop:wMbSc}(b)). However, not every combination of them is valid. In fact, it can be easily shown that if $\lim_{n\to\IF} T_n^{\beta-H} \sqrt{2\log n}=q_1\in(0,\IF)$ holds then $\lim_{n\to\IF} \frac{T_n^{H-H_0}}{ \sqrt{2\log n}}=0$. Thus, we should not expect a result like in (b.ii) and (b.iii) involving $q_1$ and $\hc_j, 2\le j\le k$ in this situation. All other combinations may be possible, as illustrated in the next example.

\end{remark}

\begin{example} \label{example:2}
As in Example \ref{example:1} this example illustrates  Theorem \ref{Thm:main2} for $T_n=(\lambda \sqln)^\gamma$, with $\gamma\in\R$ and $\lambda>0$. 
We choose to work with two representative cases. The first one is the same as in Example \ref{example:1} where $\beta>H_0>2H$. In this case, 
$$
-\IF<-\frac1{H}<-\frac{1}{H_0-H}<-\frac{1}{\beta-H}<\frac{1}{\beta-H}<\IF.
$$ 
We can show that the same eight convergence results are valid as those for $M_n$  in Example \ref{example:1}. For the second case, we consider $H_0<H<\beta<2H-H_0$. In this case,
$$
-\IF<-\frac{1}{\beta-H}<-\frac1{H}<\frac{1}{H-H_0}<\frac{1}{\beta-H}<\IF.
$$ 
A direct appliation of Theorem \ref{Thm:main2} yeilds the following convergence results, as $n\to\IF$:
\begin{itemize}
\item[(1).] If $\gamma\in\LT(-\IF, -\frac1{H}\RT)$, then $M_n \overset{d}\to 0$.
\item[(2).] If $\gamma=-\frac1{H}$, then $M_n \overset{d}\to \lambda^{-1}$.
\item[(3).] If $\gamma\in\LT( -\frac1{H}, \frac{1}{H-H_0}\RT)$, then $\sigma_0^{-1}T_n^{-H_0}(M_n - b_n) \overset{d}\to \NN$.
\item[(4).] If $\gamma=\frac{1}{H-H_0}$, then $a_n^{-1} (M_n - b_n) \overset{d}\to \Lambda+\log p_1+  \sigma_0 \lambda^{-1} \NN.$
\item[(5).] If $\gamma\in\LT( \frac{1}{H-H_0}, \frac{1}{\beta-H}\RT)$, then   $a_n^{-1} (M_n - b_n) \overset{d}\to \Lambda+\log p_1$.
\item[(6).] If $\gamma= \frac{1}{\beta-H}$ and $\lambda\in(0,\wt_0^\beta)$, then $a_n^{-1} (M_n - b_n) \overset{d}\to \Lambda+\log p_1$.
\item[(7).] If $\gamma= \frac{1}{\beta-H}$ and $\lambda=\wt_0^\beta$, then $e_n^{-1} (M_n - d_n(0)) \overset{d}\to \Lambda+\log p_1$.
\item[(8).] If $\gamma= \frac{1}{\beta-H}$ and $\lambda\in(\wt_0^\beta,\IF)$, or $\gamma> \frac{1}{\beta-H}$,  then $e_n^{-1} (M_n - d_n(\IF)) \overset{d}\to \Lambda+\log p_1$.
\end{itemize}

\end{example}

\section{Further results and proofs}

\subsection{ Proof of Proposition \ref{propTn}:} We first discuss (i). It follows, by the self-similarity, that
\BQNY
\psi_{T_u}(u)&=&\pk{\sup_{t\in[0,1]} X_H(t) -c T_u^{\beta-H} t^\beta >u T_u^{-H}   }\\
&=&\pk{\sup_{t\in[0,1]}\frac{ X_H(t)}{1+ c  t^\beta T_u^\beta/u} >u T_u^{-H}   }\\
&=&\pk{\sup_{t\in[0,1]} \frac{ X_H(t)}{1+ c   T_u^\beta/u- c  (1-t^\beta) T_u^\beta/u} >u T_u^{-H}  }\\
&=&\pk{\sup_{t\in[0,1]} \frac{ X_H(t)}{1- \frac{c   T_u^\beta/u}{1+ cT_u^\beta/u}  (1-t^\beta) } >\frac{u+cT_u^\beta} {T_u^{H}}  }.
\EQNY
Note that
\BQNY
\lim_{u\to\IF} \frac{c   T_u^\beta/u}{1+ cT_u^\beta/u} =\frac{cs_0^\beta}{1+cs_0^\beta}=c_0\ge 0.
\EQNY
Now we discuss the case $s_0>0$, implying $c_0>0$. We can easily see that, for any small $\vn>0$,  the concerned quantity $\psi_{T_u}(u)$ lies between probaiblities
\BQNY
\pk{\sup_{t\in[0,1]} Z_{\pm \vn}(t)  >\frac{u+cT_u^\beta} {T_u^{H}}  }
\EQNY
for all large enough $u$, where
\BQNY
Z_{\pm \vn}(t) = \frac{ X_H(t)}{1- c_0(1\pm \vn)  (1-t^\beta) }, \ \ \ \ t\ge0.
\EQNY
Since $s_0<t_0$,
we have that, for any sufficiently small $\vn>0,$ the variance function
\BQNY
\sigma^2_{Z_{\pm \vn}}(t) :=\text{Var}(Z_{\pm \vn}(t)) =\frac{t^{2H}}{(1-c_0(1\pm \vn)(1-t^\beta))^2},\ \ \ \ t\ge0,
\EQNY
attains its  maximum at  the unique point which is 1, and 
$\sigma_{Z_{\pm \vn}}(1)=1 $. Further,
\BQNY
1-\sigma_{Z_{\pm \vn}}(t) &=&1-\frac{t^{H}}{1-c_0(1\pm \vn)(1-t^\beta)}\\
&=&\frac{1-c_0(1\pm \vn)(1-t^\beta)-t^{H}}{1-c_0(1\pm \vn)(1-t^\beta)}\\
&=&(H-c_0\beta (1\pm \vn)) (1-t)\oo, \ \ \ \ t\uparrow 1.
\EQNY
Moreover, for the correlation function $r_{Z_{\pm \vn}}(s,t)$ of $Z_{\pm \vn}$, we have from \eqref{eq:locstaXH} that
\BQNY
1-r_{Z_{\pm \vn}}(s,t)=\frac{ 1}{2} K^2(\abs{t-s})\oo,\ \ \ \ s,t\uparrow 1.
\EQNY
Noting that
\BQNY
\lim_{u\to\IF} \frac{u+cT_u^\beta} {T_u^{H}}=\lim_{u\to\IF} \frac{u^{H/\beta} }{T_u^{H}} \frac{u} {u^{H/\beta}} \LT(1+c\frac{T_u^\beta}{u}\RT)=\IF,
\EQNY
we can apply Theorem 2.1 of \cite{DHL17b}, where the $\gamma$ defined therein is given by
\BQNY
\gamma=\lim_{t\downarrow 0} 2(H-c_0\beta (1\pm \vn)) \frac{t}{K^2(t)} = 2(H-c_0\beta(1\pm \vn)) Q.
\EQNY
Consequently, the claim in (i) for $s_0>0$ follows from an application of Theorem 2.1 in \cite{DHL17b} and  letting $\vn\to0$. 
 Next, we discuss the case  $s_0=0$, for which $c_0=0$. In this case, the concerned quantity $\psi_{T_u}(u)$ satisfies
\BQNY
\pk{\sup_{t\in[0,1]} X_{H}(t)  >\frac{u+cT_u^\beta} {T_u^{H}}  }\le \psi_{T_u}(u)\le  \pk{\sup_{t\in[0,1]}\frac{ X_{H}(t) }{1-\vn (1-t)} >\frac{u+cT_u^\beta} {T_u^{H}}  }
\EQNY
for all large enough $u$, with any small $\vn>0$. Therefore, it can be easily checked that the claim in (i) for $s_0=0$ follows similarly as the case $s_0>0.$

Now we consider (ii). The claim in (ii) follows by applying similar arguments as for Theorem 2 in \cite{HP08}. We refer to the proof of Theorem 4.1 of \cite{DHJ15} where if we set $S_v=0$ therein we  obtain the case discussed here for $x\in\R$. If $x=\IF$, 
we have, for any large $M>0$,
$$
T_u-t_0 u^{1/\beta}> M A^{1/2} B^{-1/2} u^{H/\beta+1/\beta-1}
$$
holds for all large enough $u$. Thus, applying the result with $x=M<\IF$,
\BQNY
 \Phi(M) \le \liminf_{u\to\IF}\frac{ \psi_{T_u}(u)}{ \psi_{\IF} (u)}  \le \limsup_{u\to\IF}\frac{ \psi_{T_u}(u)}{ \psi_{\IF} (u)} \le 1.
\EQNY
Letting $M\to\IF$, we obtain the required asymptotics for the case $x=\IF$.
This completes the proof.

\subsection{ Proof of Proposition \ref{prop:wMbS}} \underline{Consider (i).} For any $x\in \R$, we have
\BQNY
\pk{a_n^{-1}(\wM_n  -b_n )\le x} &=&\LT(\pk{\sup_{t\in[0,T_n]} X_i(t) -ct^\beta \le b_n +a_n x} \RT)^n\\
&=&\exp\LT( n \log\LT(1- \pk{\sup_{t\in[0,T_n]} X_i(t) -ct ^\beta> b_n +a_n x} \RT)\RT)
\EQNY
From the assumptions we in fact know for any $x\in\R$
\BQN\label{eq:bnanxT}
b_n+a_nx \sim b_n \sim (1-c\ws_0^\beta) T_n^H \sqrt{2\log n}\to\IF\ \ \ \ \  n\to\IF,
\EQN
and thus
\BQNY
\lim_{n\to\IF}\pk{\sup_{t\in[0,T_n]} X_i(t) -ct^\beta > b_n +a_n x} =0.
\EQNY
If we can show that
\BQN\label{eq:nex}
\lim_{n\to\IF} n   \pk{\sup_{t\in[0,T_n]} X_i(t) -ct > b_n +a_n x} = e^{-x},
\EQN
then
\BQNY
\pk{a_n^{-1}(\wM_n  -b_n )\le x} &\sim& \exp\LT( - n   \pk{\sup_{t\in[0,T_n]} X_i(t) -ct > b_n +a_n x}\RT)\\
&\sim&\exp\LT( - e^{-x} \RT),\ \ \ \ n\to\IF,
\EQNY
which is the required result for (i). Next, we prove \eqref{eq:nex}.
Since
\BQNY
\lim_{n\to\IF}\frac{T_n}{(b_n+a_n x)^{1/\beta}}=\LT(\frac{\ws_0^\beta}{1-c\ws_0^\beta}\RT)^{1/\beta}\in[0,t_0),
\EQNY
and 
$$
\frac{c \LT(\frac{\ws_0^\beta}{1-c\ws_0^\beta}\RT) }{1+c \LT(\frac{\ws_0^\beta}{1-c\ws_0^\beta}\RT)} =c\ws_0^\beta,
$$
we obtain, by Proposition \ref{propTn}(i), that
\BQNY
\pk{\sup_{t\in[0,T_n]} X_i(t) -ct > b_n +a_n x}  \sim f_n(b_n +a_n x) \exp\LT(-\frac{(b_n+a_n x +cT_n^\beta)^2}{2T^{2H}_n}\RT), \ \ \ \ n\to\IF,
\EQNY
where $f_n(\cdot)$ is given in \eqref{eq:tauC}.
Moreover, we have
\BQNY\frac{b_n+a_nx+cT_n^\beta}{T_n^H}&=&\sqln + \frac{1}{\sqln}\LT[x+\log \LT(f_n\LT(T_n^H\sqln-cT_n^\beta\RT)\RT)\RT]\\
&=&\sqln + \frac{1}{\sqln}\LT[ x+ \log\LT(\mathcal{D}_{c \ws_0^\beta}\LT(\sqln\RT) \LT(\sqln\RT)^{-1}\RT) \RT]\\
&\sim&\sqln,\ \ \ \ n\to\IF,
\EQNY
and, by the regular variation property of $\mathcal{D}_{c \ws_0^\beta}(\cdot)$,
\BQN\label{eq:R_nba}
f_n(b_n+a_nx)&=&\mathcal{D}_{c \ws_0^\beta} \LT( \frac{b_n+a_nx+cT_n^\beta}{T_n^H} \RT) \LT( \frac{b_n+a_nx+cT_n^\beta}{T_n^H} \RT)^{-1}\nonumber\\
&\sim & \mathcal{D}_{c \ws_0^\beta}\LT(\sqln\RT) \LT(\sqln\RT)^{-1}, \ \ \ \ n\to\IF.
\EQN
Thus, the claim in \eqref{eq:nex} follows by some elementary calculations. 

\underline{Consider (ii).} Following the same idea as in (i), we only need to show
\BQN\label{eq:nexS3}
\lim_{n\to\IF} n   \pk{\sup_{t\in[0,T_n]} X_i(t) -ct > d_n(x_0) +e_n x} = e^{-x}.
\EQN
Next, we only  focus on scenario S4, since similar arguments apply for scenario S5 as well. 
Note that, by the relations between $A, t_0$ and $\wt_0$ given in \eqref{eq:AB} and in scenario S3,
$$
\LT(\wt_0^{\beta}\sqrt{2\log n} \RT)^{\frac{1}{\beta-H}} = t_0(\Alogn)^{\frac{1}{2(\beta-H)}},
$$ and thus
\BQNY
&&\quad T_n-t_0(d_n(x_0)+e_n x)^{1/\beta}\\
&&=\LT (\wt_0^{\beta}\sqrt{2\log n}\RT)^{\frac{1}{\beta-H}} + A^{1/2} B^{-1/2}x_0(2A^2\log n)^{\frac{H+1-\beta}{2(\beta-H)}}(1+o(1)) - t_0(\Alogn)^{\frac{1}{2(\beta-H)}}\\
&&\quad\ \ -\frac{t_0A^2}{\beta-H}(\Alogn)^{\frac{2H-2\beta+1}{2(\beta-H)}}
\LT(\log\LT(R\LT(\LT(\Alogn\RT)^{1/\tau}\RT) \Phi(x_0)\RT) +x\RT)\oo\\
&&= A^{1/2} B^{-1/2}x_0(2A^2\log n)^{\frac{H+1-\beta}{2(\beta-H)}}(1+o(1)) \\
&&\ \ \ \ \ \ -
\frac{t_0A^2}{\beta-H}(\Alogn)^{\frac{2H-2\beta+1}{2(\beta-H)}}
\LT(\log\LT(R\LT(\LT(\Alogn\RT)^{1/\tau}\RT) \Phi(x_0)\RT) +x\RT)\oo,
\EQNY
as $n\to\IF.$
This implies that, under the assumption of scenario S4, 
\BQNY\label{eq:Tnt0}
\lim_{n\to\IF}\frac{T_n - t_0(d_n(x_0)+e_n x)^{1/\beta}}{A^{1/2} B^{-1/2}(d_n(x_0)+e_n x)^{H/\beta+1/\beta-1}} =x_0,
\EQNY
which also holds, with $x_0=\IF$, under the assumption of scenario S5.
Hence, by  Proposition \ref{propTn}(ii),
\BQNY
\pk{\sup_{t\in[0,T_n]} X_i(t) -ct > d_n(x_0) +e_n x} 
\sim \Phi(x_0)R(d_n(x_0) +e_n x)\exp\LT(-\frac{(d_n(x_0) +e_n x)^{\tau}}{2A^2}\RT), \ \ \ \ n\to\IF,
\EQNY
with $R(\cdot)$ given in \eqref{Ru}.
Therefore, the claim in \eqref{eq:nexS3} follows by the regular variation property of $R(\cdot)$ and using some elementary calculations  as follows
\BQNY
\frac{(d_n(x_0)+e_n x)^{\tau}}{2A^2}&=& \log n \LT( 1+ \frac{x+\log\LT(R\LT((\Alogn)^{1/\tau}\RT)\Phi(x_0)\RT)}{\tau\log n}\RT)^\tau\\
&\sim&\log n + x+\log\LT(R\LT((\Alogn)^{1/\tau}\RT)\Phi(x_0)\RT),\ \ \ n\to\IF.
\EQNY
The proof is compete.

\subsection{Proof of Theorem \ref{Thm:main1}}

\underline{Consider (a).} 
 Note that
\BQNY
\wM_n -\sup_{t\in [0,T_n]} (-\sigma_0 X(t)) \le M_n \le \wM_n +\sup_{t\in [0,T_n]} \sigma_0 X(t),
\EQNY
and, since $\lim_{n\to\IF}T_n=0$,
$$
\lim_{n\to\IF}\sup_{t\in [0,T_n]} (-\sigma_0X(t)) = \lim_{n\to\IF}\sup_{t\in [0,T_n]}\sigma_0 X(t)=0, \ \ \ \ a.s.
$$
Thus, it is sufficient to show that
\BQN\label{eq:wMk0}
\wM_n \overset{d}\to \kappa_0.
\EQN
Notice that, by an application of Theorem 2.1 in \cite{DHL17b}, 
\BQN\label{eq:XF}
\pk{\sup_{t\in[0,1]} X_i(t) >u}\  = \ \mathcal{D}_0 (u) e^{-\frac{u^2}{2}}\oo, \ \ \ u\to\IF,
\EQN
where $\mathcal{D}_0 (\cdot) $ is the regularly varying function given in \eqref{eq:DH}. Defining
\BQNY\label{eq:bnanS0}
\mu_n := \sqrt{2\log n} +\frac{\log(\mathcal{D}_0 (\sqrt{2\log n} ))}{\sqrt{2\log n}  },  \  \ \ \ \ n\in \N,   \nonumber
\EQNY
we obtain, from \eqref{eq:XF} and Proposition 2.2 in \cite{JP22}, that
\BQN \label{eq: MSXF}
\sqrt{2\log n}    \LT( \max_{i\le n}\sup_{t\in[0,1]} X_i(t)  -\mu_n \RT) \overset{d} \to \Lambda,\ \ \ \ n\to\IF.
\EQN
Next, we have
\BQNY
\max_{i\le n}\sup_{t\in[0,T_n]} X_i(t) -cT_n^\beta \le \wM_n\le \max_{i\le n}\sup_{t\in[0,T_n]} X_i(t),
\EQNY
and by self-similarity
\BQNY
\max_{i\le n}\sup_{t\in[0,T_n]} X_i(t) &\overset{d}=& T_n ^H \max_{i\le n}\sup_{t\in[0,1]} X_i(t) \\
&=&
\frac{T_n ^H}{\sqrt{2\log n} } \LT( \sqrt{2\log n}  \LT(\max_{i\le n}\sup_{t\in[0,1]} X_i(t)-\mu_n\RT)   \RT) + T_n ^H \mu_n\\
&\overset{d}\to & \kappa_0, \ \ \ \ n\to\IF.
\EQNY
Therefore, the claim in \eqref{eq:wMk0} is established and thus the proof  for (a) is complete.

Before we give the proof for (b) and (c), we shall derive some preliminary results presented in the following lemma. 


\BEL \label{Lem:Sa2} We have, 
\Ji{for any small enough $\vn_0\in(0,1)$ and any $x\in\R$,}
\begin{itemize}
\item[(i).] Under  S2 and S3,  
\BQN\label{negpartSa1}
\quad \lim_{n\to\IF}\pk{  a_n^{-1} \LT(\max_{i\le n} \sup_{0\le t\le (1-\vn_0)T_n } X_i (t) -c t^\beta  -b_n\RT) >x }=0,
\EQN
and 
\BQN\label{eq:mainInt}
a_n^{-1} \LT( \max_{i\le n}  \sup_{ (1-\vn_0)T_n \le t\le T_n} X_i (t) -c t^\beta -b_n\RT)  \overset{d} \to \    \Lambda, \ \ n\to\IF. 
\EQN

\item[(ii).] Under  S4, 
\BQN\label{negpartSa3}
\quad \lim_{n\to\IF}\pk{  e_n^{-1} \LT(\max_{i\le n} \sup_{0\le t\le (1-\vn_0)T_n } X_i (t) -c t^\beta  -d_n(x_0)\RT) >x }=0,
\EQN
and 
\BQN\label{eq:mainInt3}
e_n^{-1} \LT( \max_{i\le n}  \sup_{ (1-\vn_0)T_n \le t\le T_n} X_i (t) -c t^\beta -d_n(x_0)\RT)  \overset{d} \to \    \Lambda, \ \ n\to\IF. 
\EQN
\end{itemize}

\EEL

\prooflem{Lem:Sa2} \underline{Consider (i).}
First, note that 
\BQNY
\pk{  a_n^{-1} \LT(\max_{i\le n} \sup_{0\le t\le (1-\vn_0)T_n } X_i (t) -c t^\beta  -b_n\RT) >x } \le   n  \pk{  \sup_{0\le t\le  (1-\vn_0)T_n } X_i (t)   -c t^\beta   >  b_n+x a_n}.
\EQNY
Similarly to the proof of Proposition \ref{prop:wMbS}, we have, by Proposition \ref{propTn}(i), that
\BQNY
&& \quad   \pk{  \sup_{0\le t\le  (1-\vn_0)T_n } X_i (t)   -c t ^\beta  >  b_n+x a_n  }=  \mathcal{D}_{\wc_0}\left(\frac{b_n+x a_n+cT_n^\beta(1-\vn_0)^\beta}{T_n^H(1-\vn_0)^H}\right) \left(\frac{b_n+x a_n+cT_n^\beta(1-\vn_0)^\beta}{T_n^H(1-\vn_0)^H}\right)^{-1}\\
&& \qquad\qquad\qquad\qquad\qquad\qquad\qquad \qquad\qquad\qquad\times \exp\left(-\frac{(b_n+xa_n +cT_n^\beta(1-\vn_0)^\beta)^2}{2T_n^{2H}(1-\vn_0)^{2H}}\right)(1+o(1)),
\EQNY
as $n\to\IF,$ where
\BQN\label{eq:wc0}
\wc_0 =\frac{c (1-\vn_0)^\beta \LT(\frac{\ws_0^\beta}{1-c\ws_0^\beta}\RT) }{1+c(1-\vn_0)^\beta \LT(\frac{\ws_0^\beta}{1-c\ws_0^\beta}\RT)} =\frac{c (1-\vn_0)^\beta \ws_0^\beta }{1-c\ws_0^\beta+c(1-\vn_0)^\beta \ws_0^\beta}\ge0.
\EQN
Further, we have
\BQNY
 \frac{(b_n+xa_n +cT_n^\beta(1-\vn_0)^\beta)^2}{2T_n^{2H}(1-\vn_0)^{2H}}\sim \LT(\frac{1-c\ws_0^\beta +c\ws_0^\beta (1-\vn_0)^\beta}{(1-\vn_0)^{H}}\RT)^{2} \log n, \ \ \ \ n\to\IF,
\EQNY
and, since $c\ws_0^\beta<c\wt_0^\beta=H/\beta$, we have for all sufficiently small $\vn_0>0,$
$$
\frac{1-c\ws_0^\beta +c\ws_0^\beta (1-\vn_0)^\beta}{(1-\vn_0)^{H}} >1.
$$
These, together with the regular variation of  $\mathcal{D}_{\wc_0}(\cdot)$, imply that
\BQNY
\lim_{n\to\IF} n  \pk{  \sup_{0\le t\le  (1-\vn_0)T_n } X_i (t)   -c t^\beta   >  b_n+x a_n} =0.
\EQNY
Thus, \eqref{negpartSa1} is established. 
Next, we prove \eqref{eq:mainInt}. It is sufficient to show that, for any $x\in \R,$
\BQNY
 \lim_{n\to\IF}\pk{a_n^{-1} \LT( \max_{i\le n}  \sup_{ (1-\vn_0)T_n \le t\le T_n} X_i (t) -c t^\beta -b_n\RT) >x}=
\pk{\Lambda>x}.
\EQNY

 We have, for any $x\in \R,$ by Proposition \ref{prop:wMbS}(i),
\BQNY
\pk{\Lambda>x} &=&\lim_{n\to\IF}\pk{a_n^{-1} \LT( \wM_n -b_n\RT) >x}\\
&\ge & \limsup_{n\to\IF}\pk{a_n^{-1} \LT( \max_{i\le n}  \sup_{ (1-\vn_0)T_n \le t\le T_n} X_i (t) -c t^\beta -b_n\RT) >x}.
\EQNY
Furthermore, we have
\BQNY
\pk{a_n^{-1} \LT( \wM_n -b_n\RT) >x}&\le &\pk{a_n^{-1} \LT( \max_{i\le n}  \sup_{0\le t\le  (1-\vn_0)T_n } X_i (t) -c t^\beta -b_n\RT) >x}\\
&&+ \pk{a_n^{-1} \LT( \max_{i\le n}  \sup_{ (1-\vn_0)T_n \le t\le T_n} X_i (t) -c t^\beta -b_n\RT) >x}.
\EQNY
We have from \eqref{negpartSa1} that the first term on the right-hand side tends to 0, as $n\to\IF$. Thus,
\BQNY
\pk{\Lambda>x} \le \liminf_{n\to\IF}\pk{a_n^{-1} \LT( \max_{i\le n}  \sup_{ (1-\vn_0)T_n \le t\le T_n} X_i (t) -c t^\beta -b_n\RT) >x}.
\EQNY
This completes the proof for (i).

\underline{Consider (ii).} First, since
$$
T_n \sim \LT(\wt_0^\beta\sqrt{2\log n}\RT)^{1/(\beta-H)}\sim t_0  (d_n(x_0))^{1/\beta}, \ \ \ \ n\to\IF,
$$
the claim in \eqref{negpartSa3} follows from simialr arguments as Lemma 4.5 (see also Remark 4.6) in \cite{JP22}. 
Next, similarly to the proof of \eqref{eq:mainInt} we can prove \eqref{eq:mainInt3} by using Proposition \ref{prop:wMbS}(ii).
This completes the proof of the lemma.  \QED

\medskip
{\bf Proof of Theorem \ref{Thm:main1} continued:}
Now, we are ready to continue the proof for Theorem \ref{Thm:main1} (b)-(c) below.
We shall consider (b.i) and (b.ii)-(b.iii), respectively, followed by some arguments for (c).

\underline{Consider (b.i):}  We need to show that, for any $x\in \R$,
\BQNY
\pk{  T_n^{-H_0} (M_n -b_n)  >x} \to \pk{\sigma_0 X(1) >x},\ \ \ \ n\to\IF.
\EQNY
 For any small enough $\vn_0\in(0,1)$, we have
\BQNY 
I_1(n,\vn_0,x)\le \pk{ T_n^{-H_0} (M_n -b_n) >x} \le I_1(n,\vn_0,x) +I_2(n,\vn_0,x),
\EQNY
where
\BQNY
I_1(n,\vn_0,x)&=&\pk{T_n^{-H_0} \LT( \max_{i\le n}  \sup_{ (1-\vn_0)T_n \le t\le T_n} X_i (t) +\sigma_0X(t) -c t ^\beta-b_n\RT) >x }\\
I_2(n,\vn_0,x)&=&\pk{ T_n^{-H_0} \LT(\max_{i\le n} \sup_{0\le t\le (1-\vn_0)T_n } X_i (t) +\sigma_0X(t) -c t ^\beta -b_n\RT)  >x }.
\EQNY

For $I_1(n,\vn_0,x)$, we have
\BQNY
&&\pk{T_n^{-H_0} \LT(\max_{i\le n}  \sup_{(1-\vn_0)T_n \le t\le T_n } X_i (t)  -c t^\beta -b_n  - \sup_{(1-\vn_0)T_n \le t\le T_n} (-\sigma_0 X(t))\RT)>x}\\
&&\le I_1(n,\vn_0,x)\\
&&  \le \pk{T_n^{-H_0} \LT(\max_{i\le n}  \sup_{ (1-\vn_0)T_n \le t\le T_n} X_i (t)  -c t^\beta -b_n  + \sup_{(1-\vn_0)T_n \le t\le T_n}\sigma_0  X(t)\RT)>x} 
\EQNY
We derive by \eqref{eq:mainInt} and using the scenario assumption that,
\BQNY
&&T_n^{-H_0}\LT (\max_{i\le n}  \sup_{ (1-\vn_0)T_n \le t\le T_n} X_i (t)  -c t^\beta   -b_n\RT) \\
&& =  a_n^{-1} \LT(\max_{i\le n}  \sup_{ (1-\vn_0)T_n \le t\le T_n} X_i (t)  -c t^\beta   -b_n \RT)   (T_n^{-H_0} a_n) \\
&& \ \overset{d}\to \ 0, \quad n\to\IF.
\EQNY
Thus, by self-similarity, the independence of the Gaussian processes and the  symmetry of normal distribution, we obtain
\BQNY
&&  \pk{\sigma_0X(1) >x}= \lim_{\vn_0\to 0 }\pk{- \sup_{ 1-\vn_0 \le t\le 1} (-\sigma_0 X(t)) >x}\\
&&\le  \lim_{\vp_0\to0}\liminf_{n\to\IF} I_1(n,\vn_0,x)\le \lim_{\vp_0\to0}\limsup_{n\to\IF} I_1(n,\vn_0,x)\\ 
&&\leq \lim_{\vn_0\to 0 }\pk{ \sup_{ 1-\vn_0 \le t\le 1} \sigma_0 X(t) >x} = \pk{\sigma_0X(1) >x}.
\EQNY
To complete the proof, it is remains to show that
\BQN\label{eq:I2nvn}
\lim_{n\to\IF} I_2(n,\vn_0,x) =0.
\EQN
Since
\BQNY
I_2(n,\vn_0,x)\le   \pk{ T_n^{-H_0} \LT(\max_{i\le n} \sup_{0\le t\le (1-\vn_0)T_n } X_i (t) -c t^\beta  -b_n\RT)+\sup_{0\le t\le 1-\vn_0 }  \sigma_0X(t)  > x }
\EQNY
and $\sup_{0\le t\le 1-\vn_0  }X(t)<\IF$ \Ji{a.s.}, 
 it is sufficient to show that, for any  $x\in\R$,
\BQN\label{eq:Jnx10S1}
\lim_{n\to\IF}\pk{ T_n^{-H_0} \LT(\max_{i\le n} \sup_{0\le t\le (1-\vn_0)T_n } X_i (t) -c t^\beta  -b_n\RT) > x } = 0.
\EQN
Next, note that $T_n^{H_0}=o(b_n)$, since otherwise, there exits a subsequence $\{n_k\}_{k\geq1}$ such that $T_{n_k}^{H_0}\geq C T_{n_k}^{H}\sqrt{2\log n_k}$ holds for some positive constant $C$. Then, $T_{n_k}$ converges to $\IF$ and thus, $T_{n_k}^{1-H/\beta}/(\sqrt{2\log n_k})^{1/\beta}\geq C T_{n_k}^{1-H_0/\beta}\to\IF$, which is a contradiction with the assumption of scenarios S2 and S3.  
For the fixed $\vn_0$ and $x$,  we have by Proposition \ref{propTn}(i) 
\BQNY
&& \quad   \pk{  \sup_{0\le t\le  (1-\vn_0)T_n } X_i (t)   -c t   >  b_n+x T_n^{H_0} }=
 \mathcal{D}_{\wc_0}\left(\frac{b_n+x T_n^{H_0}+cT_n^\beta(1-\vn_0)^\beta}{T_n^H(1-\vn_0)^H}\right) \left(\frac{b_n+x T_n^{H_0}+cT_n^\beta(1-\vn_0)^\beta}{T_n^H(1-\vn_0)^H}\right)^{-1}\\
&& \qquad\qquad\qquad\qquad\qquad\qquad\qquad \qquad\qquad\qquad\times 
\exp\left(-\frac{(b_n+xT_n^{H_0} +cT_n(1-\vn_0))^2}{2T_n^{2H}(1-\vn_0)^{2H}}\right)(1+o(1)),
\EQNY
as $n\to\IF,$
with $\wc_0$ given in \eqref{eq:wc0}. 
Similarly to the proof of  \eqref{negpartSa1}, we can establish \eqref{eq:Jnx10S1}, and thus the proof for case (b.i) is finished.


\underline{Consider (b.ii) and (b.iii):} We first introduce the following notation using the indicator function:
\BQNY
I_{\{\text{case\ (b.iii)}\}} = \left\{
\begin{array}{ll}
1, &  \mbox{if\ condition\ of\ case\ (b.iii)\ is\ satisfied,\ i.e.,\ } \lim_{n\to\IF}\frac{T_n^{H-H_0}}{\sqrt{2\log n}}=q_0\in(0,\IF),\\[0.1cm]
0, &  \mbox{otherwise }.
\end{array}
\right.
\EQNY
 We need to show that, for any $x\in \R$,
\BQNY
\lim_{n\to\IF}\pk{  a_n^{-1} (M_n -b_n)  >x} = \pk{\Lambda + \sigma_0/q_0  X(1) I_{\{\text{case\ (b.iii)}\}}>x}. 
\EQNY
We will consider asymptotic upper and lower bounds, respectively.
First,   for any sufficiently small $\vn_0\in(0,1)$, we have
\BQN\label{eq:vnTn}
\pk{ a_n^{-1} (M_n -b_n) >x} &\le& \pk{a_n^{-1} \LT( \max_{i\le n}  \sup_{ (1-\vn_0)T_n \le t\le T_n} X_i (t) +\sigma_0X(t) -c t^\beta -b_n\RT) >x }\nonumber\\
&&  +\ \pk{ a_n^{-1} \LT(\max_{i\le n} \sup_{0\le t\le (1-\vn_0)T_n } X_i (t) +\sigma_0X(t) -c t^\beta   -b_n\RT)  >x }. 
\EQN
For the second term above, we have, by self-similarity,
\BQNY
&&\quad\pk{a_n^{-1}\LT( \max_{i\le n} \sup_{ 0\le t\le(1-\vn_0)T_n } X_i (t) +\sigma_0X(t) -c t^\beta  -b_n\RT) >x}\\
&&  \le \pk{a_n^{-1} \LT(\max_{i\le n}  \sup_{0\le t\le (1-\vn_0)T_n} X_i (t) -c t^\beta  -b_n  + \sup_{0\le t\le (1-\vn_0)T_n}  \sigma_0X(t)\RT)>x}\\
&& = \pk{  a_n^{-1} \LT(\max_{i\le n}  \sup_{0\le t\le (1-\vn_0)T_n } X_i (t) -c t^\beta  -b_n\RT)   + a_n^{-1} T_n^{H_0}\sup_{ 0\le t\le 1-\vn_0 }\sigma_0  X(t) >x}.
\EQNY
Under the assumption of cases (b.ii)-(b.iii), we have $\lim_{n\to\IF}a_n^{-1} T_n^{H_0}\sup_{ 0\le t\le 1-\vn_0 } \sigma_0 X(t)<\IF$, and thus by 
an application of   Lemma \ref{Lem:Sa2}(i) we derive that the above term tends to 0 as $n\to\IF$ . For the first term on the right-hand side of \eqref{eq:vnTn}, we have
\BQNY
&&\quad\pk{a_n^{-1}\LT( \max_{i\le n} \sup_{ (1-\vn_0)T_n \le t\le T_n} X_i (t) +\sigma_0X(t) -c t^\beta  -b_n\RT) >x}\\
&&  \le \pk{  a_n^{-1} \LT(\max_{i\le n}  \sup_{ (1-\vn_0)T_n \le t\le T_n} X_i (t) -c t^\beta  -b_n\RT)   + a_n^{-1} T_n^{H_0}\sup_{ 1-\vn_0 \le t\le 1} \sigma_0 X(t) >x}.
\EQNY
Thus, by Lemma \ref{Lem:Sa2}(i), 
\BQNY
&&\quad  \limsup_{n\to\IF}\pk{a_n^{-1}\LT( \max_{i\le n} \sup_{ (1-\vn_0)T_n \le t\le T_n} X_i (t) +\sigma_0X(t) -c t^\beta  -b_n\RT) >x}\\
&&\leq \lim_{\vn_0\to 0 }\pk{ \Lambda + \sigma_0/q_0\sup_{ 1-\vn_0 \le t\le 1}  X(t) I_{\{\text{case\ (b.iii)}\}} >x}\\
&& = \pk{\Lambda + \sigma_0/q_0  X(1) I_{\{\text{case\ (b.iii)}\}}>x},
\EQNY
which gives the required upper bound. For the lower bound, we have for any small $\vn_0\in(0,1)$,
\BQNY 
&&\pk{ a_n^{-1} (M_n -b_n) >x} \ge\pk{a_n^{-1} \LT( \max_{i\le n}  \sup_{ (1-\vn_0)T_n \le t\le T_n} X_i (t) +\sigma_0X(t) -c t ^\beta -b_n\RT) >x }\nonumber\\
&&\ge \ \pk{ a_n^{-1} \LT(\max_{i\le n} \sup_{ (1-\vn_0)T_n \le t\le T_n} X_i (t)  -c t^\beta   -b_n\RT)- a_n^{-1} T_n^{H_0} \sup_{ 1-\vn_0 \le t\le 1} (-\sigma_0 X(t))  >x }. \nonumber
\EQNY
By the same derivation as above, we obtain
\BQNY 
\liminf_{n\to\IF}\pk{ a_n^{-1} (M_n -b_n) >x} \ge \pk{\Lambda + \sigma_0/q_0  X(1) I_{\{\text{case\ (b.iii)}\}}>x}.
\EQNY
Hence, the proof  for (b.ii) and (b.iii) is finished.

\underline{Consider (c).}
 The proof for scenario S4 follows from the same lines as those for (b) above, by noting that
\BQNY
\lim_{n\to\IF}e_n^{-1} T_n^{H_0} = \left\{
\begin{array}{lll}
\IF, &  \mbox{if\ } 2H-H_0<\beta,\\[0.1cm]
0, &  \mbox{if\ } 2H-H_0>\beta,\\[0.1cm]
c\beta/H, & \mbox{if\ } 2H-H_0=\beta.
\end{array}
\right.
\EQNY
Next, we prove the claim for scenario S5. 
We have, for any small $\vn_0>0$, 
$$
T_n \ge 
 (1+2\vn_0)  \LT(\wt_0^\beta \sqrt{2\log n} \RT)^{1/(\beta-H)}\ge (1+\vn_0) t_0  (d_n(\IF))^{1/\beta}
$$
holds for all large enough $n.$ Since, for all large $n,$ 
\BQNY
\underset{i\le n}{\max} \ \sup_{t\in \LT[(1-\vn_0) t_0 (d_n(\IF))^{\frac{1}{\beta}},(1+\vn_0) t_0 (d_n(\IF))^{\frac{1}{\beta}} \RT]}(X_i (t)  +\sigma_0 X(t) -c t)\le M_n \le \underset{i\le n}{\max} \ \sup_{t\ge 0}(X_i (t)  +\sigma_0 X(t) -c t),
\EQNY
the proof follows from similar arguments as Theorem 3.1 in \cite{JP22}. This finishes the proof for (c), and thus completes 
the proof of the theorem. \QED

\subsection{Proof of Proposition \ref{prop:wMbSc}:} 
\underline{Consider (a).} The proof follows similarly as that for Theorem \ref{Thm:main1}(a).

\underline{Consider (b) and (c).} 
For any $x\in \R$, we have
\BQNY
\pk{a_n^{-1}(\whM_n  -b_n )\le x} &=&\prod_{i=1}^n \pk{\sup_{t\in[0,T_n]} X_i(t) -c_it^\beta \le b_n +a_n x} \\
&=&\exp\LT( \sum_{j=1}^k m_j \log\LT(1- \pk{\sup_{t\in[0,T_n]} X_i(t) -\hc_j t ^\beta> b_n +a_n x} \RT)\RT)
\EQNY
By \eqref{eq:bnanxT} we know
\BQNY
\lim_{n\to\IF}\pk{\sup_{t\in[0,T_n]} X_i(t) -\hc_j t^\beta > b_n +a_n x} =0
\EQNY
holds uniformly in $j=1,\cdots, k$. This implies, for any small $\vn>0$, that $\pk{a_n^{-1}(\whM_n  -b_n )\le x}$ lies between 
\BQNY
\exp\LT(-(1\pm\vn) \sum_{j=1}^k m_j  \pk{\sup_{t\in[0,T_n]} X_i(t) -\hc_j t ^\beta> b_n +a_n x} \RT)
\EQNY
for all large enough $n$. Similarly to the proof of Proposition \ref{prop:wMbS}, we obtain further asymptotic bounds for $\pk{a_n^{-1}(\whM_n  -b_n )\le x}$ as follows
\BQN\label{eq:exppm}
\exp\LT(-(1\pm2\vn) \sum_{j=1}^k m_j   f_n(b_n +a_n x, \hc_j) \exp\LT(-\frac{(b_n+a_n x +\hc_j T_n^\beta)^2}{2T^{2H}_n}\RT)   \RT),
\EQN
where, with $\mathcal D_{c_0}(\cdot)$ defined in \eqref{eq:DH},
\BQN \label{eq:tauCc}
 f_n(w,\hc_j):=\mathcal{D}_{c \ws_0^\beta}\left(\frac{w+\hc_j T_n^\beta}{T^H_n}\right)\left(\frac{w+\hc_j T_n^\beta}{T^H_n}\right)^{-1},  \ \ w>0.  
\EQN
Now, we focus on the summation in the exponent  of \eqref{eq:exppm}. We write (recall $c=\hc_1$)
\BQNY
&&
 I_1(n,x) :=m_1 f_n(b_n +a_n x, c) \exp\LT(-\frac{(b_n+a_n x +cT_n^\beta)^2}{2T^{2H}_n}\RT) \\
&& I_2(n,x):= \sum_{j=2}^{k} m_j f_n(b_n +a_n x, \hc_j) \exp\LT(-\frac{(b_n+a_n x +\hc_jT_n^\beta)^2}{2T^{2H}_n}\RT) 
\EQNY
Similarly to the proof of Proposition \ref{prop:wMbS}, we can obtain
\BQNY
I_1(n,x) \sim p_1 e^{-x}>0,\ \ \ \ n\to\IF.
\EQNY

Next, we discuss $I_2(n,x)$.  
We have, uniformly in $j=2,\cdots,k$,
\BQN\label{eq:bnTnH}
\frac{b_n+a_nx+\hc_jT_n^\beta}{T_n^H}&=&\sqln +(\hc_j-c) T_n^{\beta-H}+ \frac{1}{\sqln}\LT[x+\log \LT(f_n\LT(T_n^H\sqln-cT_n^\beta,\hc_j\RT)\RT)\RT]\nonumber\\
&\sim&(1+(\hc_j-c) \ws_0^\beta)\sqln,\ \ n\to\IF,
\EQN
and, by the regular variation property of $\mathcal{D}_{c \ws_0^\beta}(\cdot)$, as $n\to\IF,$
\BQN\label{eq:R_nba}
f_n(b_n+a_nx, \hc_j)&=&\mathcal{D}_{c \ws_0^\beta} \LT( \frac{b_n+a_nx+\hc_jT_n}{T_n^H} \RT) \LT( \frac{b_n+a_nx+\hc_jT_n}{T_n^H} \RT)^{-1}\nonumber\\
&\sim & \mathcal{D}_{c \ws_0^\beta}\LT((1+(\hc_j-c) \ws_0^\beta)\sqln\RT) \LT((1+(\hc_j-c) \ws_0^\beta)\sqln\RT)^{-1}.
\EQN
Below we consider the scenarios S2 and S3, seperately, to derive asymptotics for $I_2(n,x)$, as $n\to\IF$.

First consider S2, where $\ws_0 =0$. It follows that
\BQNY
\frac{(b_n+a_nx+\hc_jT_n^\beta)^2}{2T_n^{2H}}
&=&\log n+\frac{1}{2}(\hc_j-c)^2 T_n^{2(\beta-H)} +(\hc_j-c) T_n^{\beta-H} \sqln \\
&&\ \ \ \ +\LT(1+ \frac{(\hc_j-c) T_n^{\beta-H}}{\sqln} \RT) \LT[ x+\log \LT(f_n\LT(T_n^H\sqln-cT_n^\beta,\hc_j\RT)\RT) \RT]\\
&&\ \ \ \ + \frac{1}{4\log n}  \LT[ x+ \log \LT(f_n\LT(T_n^H\sqln-cT_n^\beta,\hc_j\RT)\RT) \RT]^2. 
\EQNY
Thus, some elementary calculations yield
\BQNY
 I_2(n,x)\to \left\{
\begin{array}{lll}
(1-p_1)e^{-x}, &  \mbox{if\ } \lim_{n\to\IF} T_n^{\beta-H} \sqrt{2\log n}=0,\\[0.1cm]
\sum_{j=2}^k p_j e^{-(\hc_j-c)q_1} e^{-x}, &  \mbox{if\ } \lim_{n\to\IF} T_n^{\beta-H} \sqrt{2\log n}=q_1\in(0,\IF),\\[0.1cm]
0, & \mbox{if\ } \lim_{n\to\IF} T_n^{\beta-H} \sqrt{2\log n}=\IF,
\end{array}
\right.
\ \ n\to\IF.
\EQNY
Next consider S3, where $\ws_0 \in(0,\wt_0)$.  It can be checked that, by \eqref{eq:bnTnH} and \eqref{eq:R_nba},
\BQNY
I_2(n,x) 
 \to 0, \ \ \ \ n\to\IF.
\EQNY
Now, inserting the above asymptotics for $I_1(n,x)$ and $I_2(n,x)$ into \eqref{eq:exppm} we obtain the following 
bounds for all large $n,$ 
\BQNY
\exp\LT(-(1+3\vn)   I(x)   \RT) \le \pk{a_n^{-1}(\whM_n  -b_n )\le x} \le 
\exp\LT(-(1-3\vn)  I(x)   \RT)
\EQNY
where, under S2,
\BQNY
 I(x)= \left\{
\begin{array}{lll}
e^{-x}, &  \mbox{if\ } \lim_{n\to\IF} T_n^{\beta-H} \sqrt{2\log n}=0,\\[0.1cm]
\LT(p_1+\sum_{j=2}^k p_j e^{-(\hc_j-c)q_1}\RT) e^{-x}, &  \mbox{if\ } \lim_{n\to\IF} T_n^{\beta-H} \sqrt{2\log n}=q_1\in(0,\IF),\\[0.1cm]
p_1 e^{-x}, & \mbox{if\ } \lim_{n\to\IF} T_n^{\beta-H} \sqrt{2\log n}=\IF.
\end{array}
\right.
\EQNY
and, under S3,
\BQNY
I(x)= p_1 e^{-x}.
\EQNY
Therefore, the claims in (b) and (c) are established by letting $\vn\to0$.

\underline{Consider (d).} Since the proof under S5 is similar to the proof under S4, we shall only present a proof  under S4.  Similarly as for (c) above,
we conclude, for any small $\vn>0$, that $\pk{e_n^{-1}(\whM_n  -d_n(x_0) )\le x}$ lies between 
\BQNY
\exp\LT(-(1\pm\vn) \sum_{j=1}^k m_j  \pk{\sup_{t\in[0,T_n]} X_i(t) -\hc_j t ^\beta> d_n(x_0) +e_n x} \RT)
\EQNY
for all large enough $n$. Similarly to the proof of Proposition \ref{prop:wMbS}, by an application of Proposition \ref{propTn}(ii) we obtain asymptotic bounds for $\pk{e_n^{-1}(\whM_n  -d_n(x_0) )\le x}$ as follows
\BQNY 
\exp\LT(-(1\pm2\vn) \sum_{j=1}^k m_j \Phi(x_0) R(d_n(x_0) +e_n x)\exp\LT(-\frac{(d_n(x_0) +e_n x)^{\tau}}{2A^2}\RT)   \RT).
\EQNY
Note that $A=A(c)$ defined in \eqref{eq:AB}, as a function of $c$, is strictly dicreasing. Thus, using similar arguments as before we can estiblish the result under S4.
This completes the proof of the proposition.

\subsection{Proof of Theorem \ref{Thm:main2}}

The proof follows by the same arguments as those for Theorem \ref{Thm:main1}, for which we should use the Proposition \ref{prop:wMbSc} as a replacement of Proposition \ref{prop:wMbS}.

\bigskip

{\bf Acknowledgement}:
The research of Xiaofan Peng is partially supported by  National Natural Science Foundation of China (11701070, 71871046).

\bibliographystyle{plain}
\bibliography{gausbibruinABCD}

\bigskip

\end{document}